\documentclass[twocolumn]{article}
\usepackage{imr}
\usepackage{amsmath,amssymb,amscd}
\usepackage{graphicx}

\usepackage[
colorlinks=true,
citecolor=blue,
urlcolor=blue,
linktocpage,
pagebackref
]{hyperref}

\def\bx{\boldsymbol{x}}
\def\bxs{\boldsymbol{x_s}}
\def\bbx{\boldsymbol{\bar{x}}}

\newtheorem{remark}{Remark}

\begin{document}

\title{Adaptive Surface Fitting and Tangential Relaxation for
High-Order Mesh Optimization\footnote{Performed under the auspices of the U.S. Department
of Energy under Contract DE-AC52-07NA27344 (LLNL-CONF-819631).}}
\author{
Patrick Knupp$^1$ \and
Tzanio Kolev$^2$ \and
Ketan Mittal$^2$ \and
Vladimir Z. Tomov$^{2,\dagger}$}
\date{
$^1$Dihedral LLC, Bozeman, MT, U.S.A. knupp.patrick@gmail.com \\
$^2$Lawrence Livermore National Laboratory, Livermore, CA, U.S.A. \\
kolev1@llnl.gov, mittal3@llnl.gov, tomov2@llnl.gov \\
$^\dagger$Corresponding author}

\abstract{
We propose a new approach for controlling the characteristics of
certain mesh faces during optimization of high-order curved meshes.
The practical goals are tangential relaxation along initially aligned
curved boundaries and internal surfaces, and mesh fitting to initially
non-aligned surfaces.
The distinct feature of the method is that it utilizes discrete
finite element functions (for example level set functions) to define implicit
surfaces, which are used to adapt the positions of certain mesh nodes.
The algorithm does not require CAD descriptions or analytic parametrizations,
and can be beneficial in computations with dynamically changing geometry,
for example shape optimization and moving mesh multimaterial simulations.
The main advantage of this approach is that it completely avoids geometric
operations (e.g., surface projections), and all calculations can be
performed through finite element operations.
}

\keywords{
high-order meshes, node movement,
tangential relaxation, surface fitting, finite elements}

\maketitle
\thispagestyle{empty}
\pagestyle{empty}

%-------------------------------------------------

\section{Introduction}
High-order finite element (FE) methods are becoming increasingly important in
computational science due to their potential for better simulation accuracy and
favorable scaling on modern architectures \cite{Fischer2002,fischer2020scalability,anderson2021mfem, CEED2021}.
A vital component of such methods is the use of high-order representation
for the geometry, represented by a high-order computational mesh.
Such meshes are essential for achieving the optimal convergence rates on domains
with curved boundaries/interfaces, symmetry preservation, and alignment with
the physics flow in moving mesh simulations
\cite{Shephard2011, Dobrev2012, Boscheri2016}.

In order to fully utilize the benefits of high-order geometry representation,
however, one must be able to
control the quality and adapt the properties of a high-order mesh.
Two common requirements for mesh optimization methods are (1) to fit certain
mesh faces to a given surface representation
(see example in Section \ref{sec_fit_2D}),
and (2) to perform tangential node movement along a mesh surface
(see example in Section \ref{sec_tg}).
This paper is concerned with these two requirements, in the particular case
when the surface representation is defined by a discrete (or implicit) function
without an analytic parametrization.
Common examples of this scenario include the use of level set functions to
represent material interfaces in multimaterial simulations, or to represent
an evolving geometry in topology and shape optimization applications.

Aligning mesh faces to curved boundaries through FE-based variational
formulations is a common approach in the mesh generation literature.
A non-exhaustive list of publications and recent advances on the subject is
given by \cite{Roca2016, Toulorge2016, Roca2017, Sarrate2016, Moxey2015, Sevilla2014, Persson2016, Roca2019}.
The proposed algorithm in this paper falls into the category
of variational methods that force surface fitting incrementally
through variational penalty terms.
These methods enforce the surface fitting weakly, thus allowing more freedom
for the boundary nodes and natural tangential sliding around the surface.
Specific examples for such methods include indirect utilization of the CAD
parametrization through periodic surface projections \cite{Roca2017}, and
weak enforcement through Lagrange multipliers \cite{Roca2019}.
The common theme among all of the above approaches is the use of
analytic CAD parametrization of the surfaces.
This paper is explicitly focused on the case when the surface of interest is
known implicitly, i.e., there is no parametrization and the surface is
prescribed only through a discrete finite element function.
Other related works in the area of implicit surface fitting and
tangential relaxation include the \textit{DistMesh} algorithm that generates a
Delaunay triangulation in the domain of interest, followed by solving for force
equilibrium in the Delaunay structure to obtain a body-fitted linear mesh \cite{persson2004simple};
Rangarajan's method to generate a boundary-fitting triangulation by trimming a
conforming (low-order) mesh and projecting the boundary vertices of the trimmed
mesh to the desired level set \cite{rangarajan2019algorithm};
Chen's interface-fitted (linear) mesh generator that uses geometric operations
such as splitting and merging to modify existing elements in a mesh to align
them to an interface \cite{chen2017interface};
and Mittal's distance function-based approach for approximate tangential
relaxation for surface nodes \cite{mittal2019mesh}.

In this paper we propose a new approach that is applicable to both
tangential relaxation and surface fitting.
The mesh optimization problem is formulated as a variational minimization of a
chosen mesh-quality metric, through the Target-Matrix Optimization Paradigm
\cite{ETHOS2019, Knupp2012}, with additional penalty terms that enforce
the desired tangential motion and/or surface alignment.
These penalty terms connect the concept of mesh motion to the discrete
finite element function that defines the desired node position.
The method utilizes a single objective function for all mesh nodes, that is,
the nodes that are selected for alignment move together with all other nodes.
As the penalty terms depend on discrete functions, there is an interpolation
procedure that makes the functions available on different meshes \cite{IMR2019}.
The optimization method is based on global node movement and does not
alter the topology of the starting mesh.

The main advantage of the proposed approach is that it completely avoids
geometric operations (e.g., surface projections), and all calculations
can be performed through \emph{high-order} finite element operations.
Thus the main steps of the method are independent of dimension,
order of the mesh, and types of elements.
Another benefit is that, unlike most geometric operations, the FE-based surface
fitting terms can be differentiated, s.t. there is no need for an outside loop
around the main nonlinear solver. The drawback is that the notions of alignment
and fitting are always approximate (imposed weakly), and elimination of small
surface features could remain unnoticed by the nonlinear solver.
Thus, the method may require additional modifications in situations that require
the exact preservation of certain features, e.g.,
for representation of sharp corners in 3D.

The rest of the paper is organized as follows.
In Section \ref{sec_prelims} we review the basic TMOP components and
our framework to represent and optimize high-order meshes.
The technical details of the proposed method for surface fitting and tangential
relaxation are described in Section \ref{sec_fit}.
Section \ref{sec_examples} presents several academic tests that
demonstrate the main features of the methods,
followed by a conclusion in Section \ref{sec_concl}.

%-------------------------------------------------

\section{Preliminaries}
\label{sec_prelims}

The method presented in this paper is an extension of our previous work on
the TMOP framework for high-order meshes \cite{ETHOS2019} and its extension
to simulation-driven adaptivity \cite{ETHOS2020}.
In this section we summarize the main concepts and notation that are related to
the understanding of the newly developed algorithms.

The domain $\Omega \in \mathbb{R}^d$ is
discretized as a union of curved mesh elements of order $k$.
Discrete representation of these elements is obtained by utilizing a set of
scalar basis functions $\{ \bar{w}_i \}_{i=1}^{N_w}$ on the reference element
$\bar{E}$.
This basis spans the space of all polynomials of degree at most $k$ on the
given element type (quadrilateral, tetrahedron, etc.).
The position of an element $E$ in the mesh $\mathcal{M}$
is fully described by a matrix
$\mathbf{x}_E$ of size $d \times N_w$ whose columns represent the coordinates
of the element control points ({\em nodes} or element {\em degrees of freedom}).
Given $\mathbf{x}_E$, we introduce the map $\Phi_E:\bar{E} \to \mathbb{R}^d$
whose image is the element $E$:
\begin{equation}
\label{eq_x}
\bx(\bbx) =
   \Phi_E(\bbx) \equiv
   \sum_{i=1}^{N_w} \mathbf{x}_{E,i} \bar{w}_i(\bbx)\,,
   \qquad \bbx \in \bar{E},
\end{equation}
where we used $\mathbf{x}_{E,i}$ to denote the $i$-th column of $\mathbf{x}_E$,
i.e., the $i-$th degree of freedom of element $E$.
The Jacobian $A_{d \times d}$ of the mapping \eqref{eq_x} at any
reference point $\bar{x} \in \bar{E}$ is computed as
\begin{equation}
\label{eq_A}
  A_E(\bar{x}) = \frac{\partial \Phi_E}{\partial \bar{x}} =
    \sum_{i=1}^{N_w} \mathbf{x}_{E,i} [ \nabla \bar{w}_i(\bar{x}) ]^T \,.
\end{equation}
The $\mathbf{x}_E$ control points of every element
are arranged in a global vector $\bx$ of size $d \cdot N_x$
that stores the coordinates of all node positions.
Throughout this paper we exclusively work with global vectors and use the
following notation:
\begin{equation}
  \bx = (x_1 \dots x_d)^T, ~
  x_a = \sum_{i=1}^{N_x} x_{a,i} w_i(\bbx), ~ a = 1 \dots d.
\end{equation}

\begin{figure}[b!]
\centerline
{
  \includegraphics[width=0.4\textwidth]{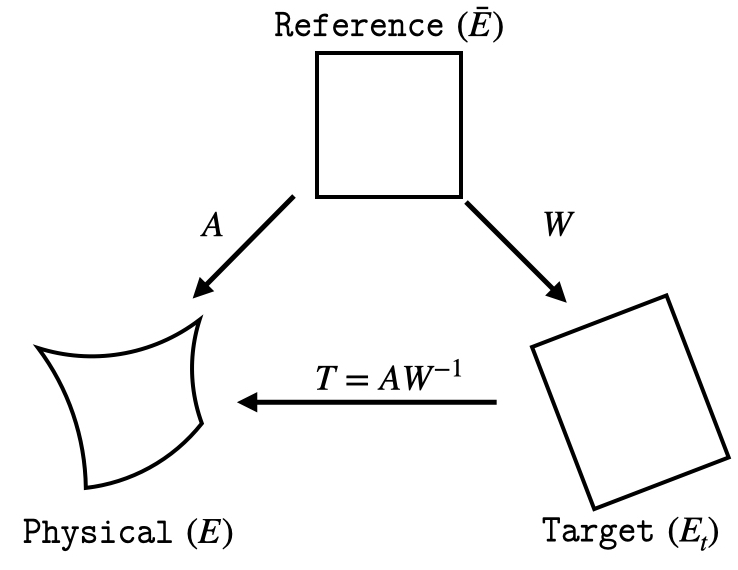}
}
\caption{Schematic representation of the target to physical transformation in TMOP.}
\label{fig_tmop}
\end{figure}

The positions of the mesh nodes are optimized by minimizing a
global objective function:
\begin{equation}
\label{eq_F}
  F(\bx) = \sum_{E \in \mathcal{M}} \int_{E_t} \mu(T(\bx)) d \bx_t +
           F_{\sigma}(\bx),
\end{equation}
where $E_t$ are user-defined target elements through their Jacobian matrices
$W_{d \times d}$;
$T = A W^{-1}$ is the weighted Jacobian matrix from
\emph{target} to \emph{physical} coordinates, see Figure \ref{fig_tmop};
$\mu(T)$ is a mesh quality metric that defines a measure of the difference
between the target and the physical geometric properties at a given location.
The above configuration is the classical TMOP setup that is used as a backbone
of many other mesh optimization methods \cite{Roca2016, Peiro2018}.
The term $F_{\sigma}(\bx)$ is the subject of this paper and will
be discussed in Section \ref{sec_fit}.
This term is used to enforce both tangential relaxation and surface fitting.

The objective function \eqref{eq_F} is minimized by solving
$\partial F (\bx) / \partial \bx = \boldsymbol{0}$.
This is done by a global nonlinear solve that modifies all values in $\bx$ at once.
Each iteration is enhanced by a line search procedure to ensure mesh validity,
as described in \cite{dobrev2021hr}. Note that we typically use the Newton's method,
but occasionally we also employ the limited memory Broyden--Fletcher--Goldfarb--Shanno
(L-BFGS) method. The convergence criterion for the solver is based on the norm of
the gradient of the objective function with respect to the current and the initial mesh,
$|\nabla F(\bx)|/|\nabla F(\bx_0)| \leq \varepsilon$, and we set $\varepsilon=10^{-6}$ for
all the numerical results in this paper.
For the case of Newton, the linear system inside each iteration is inverted
by the standard minimum residual (MINRES) algorithm
with an $l_1$ - Jacobi preconditioner.

%-------------------------------------------------

\section{Adaptive Surface Fitting and Tangential Relaxation}
\label{sec_fit}

The main input of the method is a scalar FE function
$\sigma_0(\boldsymbol{x_0})$ defined with respect to the
initial mesh $\mathcal{M}_0$.
The zero level set of this function specifies the surface of interest.
In the following discussion we focus explicitly on the discrete case, but
the method is also applicable when $\sigma$ is given analytically.
The proposed algorithm is nearly identical for the cases of
fitting and tangential relaxation; we will generally not separate the two cases,
and the differences between them will be noted explicitly.

A starting step of the algorithm is to choose a subset of the mesh nodes,
$\mathcal{S}$, which will be constrained on the zero level set.
For the case of tangential relaxation (when we try to preserve a mesh surface),
$\mathcal{S}$ is known by definition,
as the surface of interest is defined by some of the faces of $\mathcal{M}_0$.
The set $\mathcal{S}$ is also known when one wants to fit boundary faces of a
non-aligned initial mesh.
For the case of fitting a non-aligned initial mesh to an internal surface,
however, choosing the set $\mathcal{S}$ can be a nontrivial problem, which
we do not address in this paper. In our fitting tests we
choose $\mathcal{S}$ through heuristics related to the shape of interest.

\begin{figure}[b!]
\centerline
{
  \includegraphics[width=0.4\textwidth]{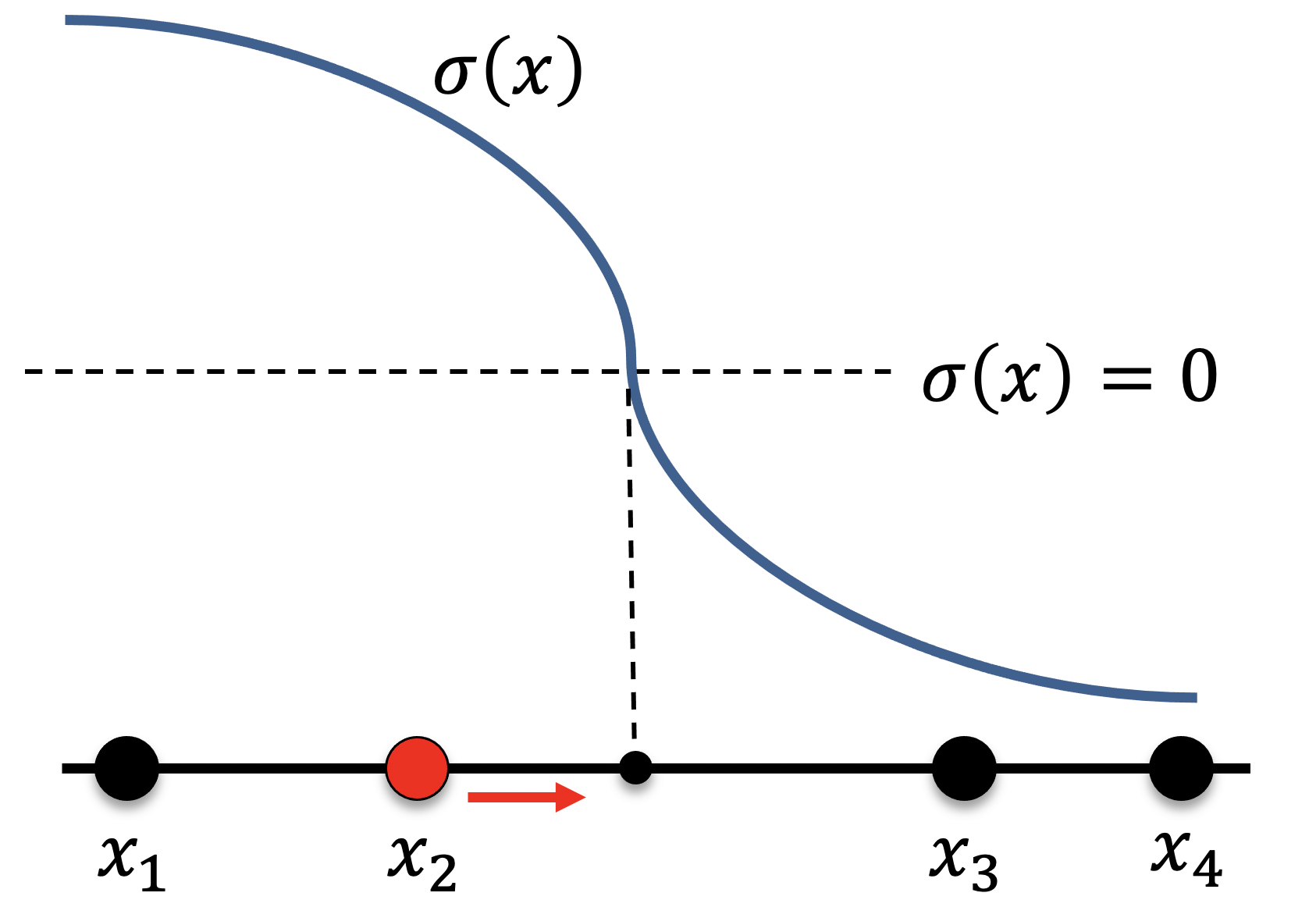}
}
\caption{Schematic representation of aligning mesh nodes to
         the zero level set of a discrete function.}
\label{fig_sigma1D}
\end{figure}

Once the set $\mathcal{S}$ is determined, the main idea of the adaptive surface
alignment is to (i) move the mesh in a manner that leads to $\sigma(\bxs)$ = 0,
for all $s \in \mathcal{S}$, i.e., to place all nodes $s \in \mathcal{S}$ as
close as possible to the zero level set, and (ii) to maintain optimal mesh
quality as defined by the quality metric $\mu$.
As there is no pre-determined unique target position for each node of
$\mathcal{S}$, the method naturally allows tangential relaxation along
the interface of interest.
Simple 1D scheme illustration of the approach is shown in
Figure \ref{fig_sigma1D}, where the mesh node
$\boldsymbol{x_2}$ is fitted to the zero level set of $\sigma$.
For the 2D and 3D cases the level set becomes a curve and a surface,
respectively.

The method is called \textit{adaptive}, as the optimal mesh positions are
obtained from a discrete function that comes from the simulation of interest.
One of the main difficulties in this case comes from the fact that
$\sigma$ is only defined with respect to the initial mesh.
For a node $s \in \mathcal{S}$, let $\boldsymbol{x_{s0}}$ be its position in the
initial mesh $\mathcal{M}_0$,  and let $\bxs$ be its corresponding mesh position
in the current mesh $\mathcal{M}$ which is evolving during the optimization.
Generally $\sigma_0(\boldsymbol{x_{s0}}) \neq 0$, as the chosen nodes may not
originally be on the zero level set.
After each mesh modification, we recompute a function $\sigma(\bx)$ on the
new mesh (see additional details in Section \ref{sec_interpolation}), so that
ideally $\sigma(\bx) = \sigma_0(\bx)$ for all points $\bx$, which represents
equivalence of $\sigma_0$ and $\sigma$ physical space.

Let the FE representation of $\sigma$ be
\[
  \sigma(\bx) = \sum_{i=1}^{N_{\phi}} \sigma_i \phi_i(\bx),
\]
where $\{\phi\}_{i=1}^{N_{\phi}}$ are the basis functions and
$N_{\phi}$ is the total number of DOFs of $\sigma$.
We make two important assumptions about the FE space of $\sigma$, namely:
\begin{itemize}
\item The FE space of $\sigma$ uses interpolatory basis functions, i.e.,
      $\phi_j(\boldsymbol{x_i})=\delta_{ij}$ for all nodes $\boldsymbol{x_i}$,
      so that $\sigma(\boldsymbol{x_i})$ = 0 if and only if $\sigma_i = 0$.
\item The FE nodes of $\sigma$ and $\bx$ coincide.
      For example, if the mesh control points are based on the Gauss-Lobatto
      reference positions (as explained in Section \ref{sec_prelims}), then the
      FE space of $\sigma$ must be based on the Gauss-Lobatto positions as well.
\end{itemize}
If the application provides a $\sigma$ function that does not satisfy the
above, a FE projection operation on an appropriate space can be performed.
The above assumptions provide a very convenient property, namely,
that minimizing the absolute value of a FE coefficient $\sigma_s$,
for a node $s \in \mathcal{S}$, would lead to approaching $\sigma(\bxs) = 0$
which would mean that the mesh node $\bxs$ is near the zero level set of
$\sigma$; in other words, $\bxs$ would be aligned to the surface of interest.
This property is used to formulate the penalty term $F_{\sigma}$ in \eqref{eq_F}.

To incorporate the idea of minimizing $\sigma_s$ into the objective function
\eqref{eq_F}, we first define the \textit{restricted} FE function
$\bar{\sigma}(x)$, s.t. for each FE coefficient $i$ we have:
\begin{equation}
\label{eq_sbar}
  \bar{\sigma}_i =
  \begin{cases}
     \sigma_i & \text{if} ~ i \in \mathcal{S}, \\
     0        & \text{otherwise}.
  \end{cases}
\end{equation}
In other words,
$\bar{\sigma}(x) = \sum_{s \in \mathcal{S}} \sigma_s \phi_s(x)$.
The extra penalty term in the objective function \eqref{eq_F} is formulated as:
\begin{equation}
\label{eq_F_sigma}
  F_{\sigma} =
    \frac{\omega_{\sigma}}{c_{\sigma}}
    \int_{E_t} \bar{\sigma}(x)^2,
\end{equation}
where $\omega_{\sigma}$ is a user-defined weight,
and $c_{\sigma}$ is a normalization constant.
The normalization scaling is used to ensure that $F_{\sigma}$ is invariant
with respect to mesh refinement and scaling of the domain:
\[
  \frac{1}{c_{\sigma}} =
  \begin{cases}
     \frac{1}{N_E} \frac{1}{V_{\Omega, avg}} = \frac{1}{V_{\Omega}}
       & \text{for volumetric targets}, \\
     \frac{1}{N_E}
       & \text{otherwise},
  \end{cases}
\]
where $N_E$ and $V_{\Omega}$ are the total number of mesh elements and volume
of the domain, respectively. Here by \textit{volumetric targets} we refer to
the target Jacobian matrices $W$ that contain volumetric information and
their determinants have unit of volume.

The extra term in the objective function penalizes the nonzero
values of $\sigma(x_s)$ for all $s \in \mathcal{S}$.
Minimizing this term represents weak enforcement of $\sigma(\bxs) = 0$,
only for the nodes in $\mathcal{S}$, while ignoring the values of $\sigma$ for
the nodes outside $\mathcal{S}$.
Note that minimizing the final nonlinear objective function
$F = F_{\mu} + F_{\sigma}$, see \eqref{eq_F}, would treat all nodes together,
i.e., the nonlinear solver would make no explicit separation between surface
and non-surface nodes.

\begin{remark}
Since the FE basis of $\bx$ is also interpolatory, we
can simply assume that the basis functions of $\sigma$ and $\bx$ are the same,
that is, $\phi \equiv w$. Throughout the paper we keep both notations to
distinguish the terms related to $\sigma$ from the ones related to the mesh
positions.
\end{remark}

%-------

\subsection{Derivatives}

As our default choice for nonlinear optimization is the Newton's method, we
must compute two derivatives of $F_{\sigma}$ with respect to the mesh nodes.
In this discussion we're exclusively targeting the discrete case, when the level
set function $\sigma(\bx)$ is a discrete FE function.
All derived formulas are also applicable to the analytic case; a remark about
the required modification is given at the end of the section.

Let the FE expansion for $\sigma$ be
$\sigma(\bx) = \sum_k \sigma_k \phi_k(\bx)$,
and let the FE position function be $\bx = (x_1 \dots x_d)^T$ where $d$ is the dimension and each component can be written as $x_a = \sum_i x_{a,i} w_i(\bbx)$.
For a node $s \in \mathcal{S}$, let $\bxs$ be its position.
As the chosen FE basis is interpolatory, we have
$\bxs = (x_{1,s} \dots x_{d,s})^T$, and we can write the following equivalence:
\[
\sigma_s = \sigma(\bxs) = \sum_{k=1}^{N_x} \sigma_k \phi_k(\bxs),
  \quad \text{since} \quad
\phi_k(\bxs) = \delta_{ks}.
\]
Thus we can rewrite $\bar{\sigma}$ as
\[
\bar{\sigma}(\bx) =
  \sum_{s \in \mathcal{S}} \sigma_s \phi_s(\bx) =
  \sum_{s \in \mathcal{S}} \sum_{k=1}^{N_x} \sigma_k \phi_k(\bxs) \phi_s(\bx).
\]
The above formula is useful for derivative computations.
Namely, we have $\phi_k(\bxs) = 0$ when $k \neq s$, however, the derivatives
$\partial \phi_k(\bxs) / \partial \bx \neq 0$.
A chain rule for the above expression (given below) expresses how the change of
the DOF values $\sigma(\bxs)$ contribute to the change of the quadrature point
values of the restricted FE function $\bar{\sigma}(\bx)$.
The formulas for the first and second derivatives of $F_{\sigma}$
are the following:
\[
\begin{split}
\frac{\partial F_{\sigma}}{\partial x_{a,i}}
= \frac{2 \omega_{\sigma}}{c_{\sigma}} \int_{E_t} &
  \bar{\sigma}(\bx)
  \frac{\partial \bar{\sigma}(\bx)}{\partial x_a}
  \frac{\partial x_a}{\partial x_{a,i}} \\
= \frac{2 \omega_{\sigma}}{c_{\sigma}} \int_{E_t} &
  \bar{\sigma}(\bx)
  \sum_{s \in \mathcal{S}} \sum_k \sigma_k \bigg(
    \frac{\partial \phi_k(\bxs)}{\partial x_a} \phi_s(\bx) + \\
    & \qquad \qquad \qquad
    \phi_k(\bxs) \frac{\partial \phi_s(\bx)}{\partial x_a} \bigg)
  w_i(\bbx),
\end{split}
\]
\[
\begin{split}
\frac{\partial^2 F_{\sigma}}{\partial x_{b,j} \partial x_{a,i}}
= \frac{2 \omega_{\sigma}}{c_{\sigma}} \int_{E_t} &
  \bigg(
    \frac{\partial \bar{\sigma}(\bx)}{\partial x_b}
    \frac{\partial \bar{\sigma}(\bx)}{\partial x_a} + \\
    & \quad
    \bar{\sigma}(\bx)
    \frac{\partial^2 \bar{\sigma}(\bx)}{\partial x_b \partial x_a}
  \bigg)
  \frac{\partial x_a}{\partial x_{a,i}} \frac{\partial x_b}{\partial x_{b,i}} \\
= \frac{2 \omega_{\sigma}}{c_{\sigma}} \int_{E_t} &
  \left( \mathcal{D}_a \mathcal{D}_b +
         \bar{\sigma}(\bx) \mathcal{D}^2
  \right)
  w_i(\bbx) w_j(\bbx),
\end{split}
\]
where
\[
\begin{split}
\mathcal{D}_* =
  \sum_{s \in \mathcal{S}} \sum_k & \sigma_k \left(
    \frac{\partial \phi_k(\bxs)}{\partial x_*} \phi_s(\bx) +
    \phi_k(x_s) \frac{\partial \phi_s(\bx)}{\partial x_*} \right), \\
\mathcal{D}^2 =
  \sum_{s \in \mathcal{S}} \sum_k & \sigma_k
  \bigg( \frac{\partial \phi_k(\bxs)}{\partial x_a}
         \frac{\partial \phi_s(\bx)}{\partial x_b} +
         \frac{\partial^2 \phi_k(\bxs)}{\partial x_b \partial x_a} \phi_s(\bx)~+ \\
         & \quad
         \frac{\partial \phi_k(\bxs)}{\partial x_b}
         \frac{\partial \phi_s(\bx)}{\partial x_a} +
         \phi_k(\bxs)
         \frac{\partial^2 \phi_s(\bx)}{\partial x_b \partial x_a}
  \bigg),
\end{split}
\]
\[
a,b = 1 \dots d, \quad i,j = 1 \dots N_x.
\]
The above integrals are approximated by a standard Gauss-Lobatto integration
rule of order that depends on the used mesh degree.
Note that the above formulas mix gradients at quadrature points,
e.g., $\partial \phi_s(\bx) / \partial x_a$, and gradients at the nodes of
$\mathcal{S}$, e.g., $\partial \phi_k(\bxs) / \partial x_a$.
Furthermore, note that gradients of the basis functions are with
respect to physical coordinates, s.t. integration in reference space would
require to compute $\nabla \phi(\bx) = A^{-T} \nabla \hat{\phi}(\bbx)$.
This leads to more involved computations of the second derivatives.
One possibility to avoid this is to utilize approximate second derivatives,
which are obtained through repeated application of the FE discrete gradient
operator.

\noindent
{\bf Analytic case.}
When $\sigma(\bx)$ is prescribed analytically, the above formulas still hold
with a slight simplification.
In this case it is more convenient to rewrite
the restricted function $\bar{\sigma}(\bx)$ as
\[
\bar{\sigma}(\bx) =
  \sum_{s \in \mathcal{S}} \sigma(\bxs) \phi_s(\bx).
\]
Then the first derivative of $F_{\sigma}$ becomes:
\[
\begin{split}
\frac{\partial F_{\sigma}}{\partial x_{a,i}}
= \frac{2 \omega_{\sigma}}{c_{\sigma}} \int_{E_t} &
  \bar{\sigma}(\bx)
  \frac{\partial \bar{\sigma}(\bx)}{\partial x_a}
  \frac{\partial x_a}{\partial x_{a,i}} \\
= \frac{2 \omega_{\sigma}}{c_{\sigma}} \int_{E_t} &
  \bar{\sigma}(\bx)
  \sum_{s \in \mathcal{S}} \bigg(
    \frac{\partial \sigma(\bxs)}{\partial x_a} \phi_s(\bx) + \\
    & \qquad \qquad
    \sigma(\bxs) \frac{\partial \phi_s(\bx)}{\partial x_a} \bigg)
  w_i(\bbx),
\end{split}
\]
where all derivatives $\partial \sigma(\bx)/\partial x_a, a = 1 \dots d$
are known analytically.
The formulas for the second derivatives are modified in a similar manner.

%-------

\subsection{Interpolation Between Meshes}
\label{sec_interpolation}

As we always use an iterative method that produces a sequence of meshes, we
need to compute $\sigma(\bx)$ and its gradients on any of these meshes.
This is straightforward when $\sigma$ is given analytically, but in practice,
$\sigma$ is often a discrete function produced by a numerical simulation.
As such, it is only defined with respect to the initial mesh $\mathcal{M}_0$.
In these cases the ability to compute $\sigma(\bx)$ on an evolved mesh
$\mathcal{M}$ is a major requirement for the algorithm.
An exact equivalence in physical space would mean to have
$\sigma(\bx) = \sigma_0(\bx)$ for all points $\bx$, but this is not possible
due to the discrete FE representation of these functions.

The reconstruction of $\sigma$ on $\mathcal{M}$ can be performed entirely by
FE operations.
Using the topological equivalence of $\mathcal{M}_0$ and $\mathcal{M}$, one can
define mesh velocity $\boldsymbol{v} = \bx - \boldsymbol{x_0}$ and solve the following advection PDE in pseudo-time $\tau \in [0,1]$:
\[
  \frac{d \sigma(\boldsymbol{x}_{\tau}, \tau)}{d \tau} =
    \boldsymbol{u} \nabla \sigma(\boldsymbol{x}_{\tau}, \tau) \cdot , \quad
  \sigma(\boldsymbol{x}_0, 0) = \sigma_0(\boldsymbol{x}_0),
\]
where $\boldsymbol{x}_{\tau} = \boldsymbol{x}_0 + \tau \boldsymbol{v}$.
Further details about this procedure can be found in Section 4.2 of
\cite{IMR2019}.

Another option is high-order interpolation between meshes in physical space,
which is enabled by the open-source library, \textit{gslib} \cite{gslibrepo}.
The \emph{findpts} set of routines in \emph{gslib} provide two key functionalities.
First, for a given set of points in physical space, \emph{findpts} determines
the computational coordinates for each point, i.e.,
the element $E$ which contains the point and the reference-space
coordinate $\bbx$ inside the corresponding reference element $\bar{E}$.
Second, using the computational coordinates, any given high-order FE function is interpolated \eqref{eq_x}.
These two key functionalities allow the transfer of a high-order FE function from
one mesh onto another.
The reader is referred to Section 2.3 of \cite{Mittal2019} for further details.

%-------------------------------------------------

\section{Numerical Examples}
\label{sec_examples}

In this section we demonstrate the main properties of the method on
several proof-of-concept 2D and 3D tests.
We start with standard surface fitting tests on different element types,
followed by treatment of smooth and non-smooth internal interfaces.
Extensive evaluation of the proposed method on more complicated geometries
and practical problems will be performed in future work, as these cases
require additional capabilities that are not established yet in our software.
Nevertheless, the proposed method can be readily utilized
as a building block by other established mesh optimization frameworks like
\cite{Roca2016, Peiro2018}.

Unless notes otherwise, the presented tests utilize the following composite
mesh quality metrics:
\begin{equation}
\begin{split}
  \mu_{80}  &= (1-\gamma) \mu_2     + \gamma \mu_{77}, \\
  \mu_{333} &= (1-\gamma) \mu_{302} + \gamma \mu_{316},
\end{split}
\end{equation}
where $\gamma = 0.5$ and
\[
\begin{split}
  \mu_2     &= 0.5 \frac{|T|^2}{\tau} - 1, \quad \quad
  \mu_{77}   = 0.5 \left(\tau - \frac{1}{\tau} \right)^2, \\
  \mu_{302} &= \frac{|T|^2|T^{-1}|^2}{9} - 1, \quad
  \mu_{316}  = 0.5 \left( \tau + \frac{1}{\tau} \right) - 1,
\end{split}
\]
where $|T|$ and $\tau$ are the Frobenius norm and determinant of $T$, respectively.
The metric $\mu_{80}$ is a 2D \textit{shape+size} metric, while
$\mu_{333}$ is a 3D \textit{shape+size} metric.
Both are polyconvex in the sense of \cite{Garanzha2010, Garanzha2014},
i.e., the metric integral $F_{\mu}$ in \eqref{eq_F} theoretically has a minimizer.
A thorough investigation of the theoretical properties
of the above (and many other) mesh quality metrics and metric types
can be found in \cite{Knupp2020}.
Exploring how the smoothness properties of $\sigma$ affect the convexity
properties of the full objective function $F = F_{\mu} + F_{\sigma}$
will be the subject of future studies.

Our implementation utilizes the MFEM finite element library \cite{anderson2021mfem}.
This implementation is freely available at \url{https://mfem.org}.

%-------

\subsection{2D Surface Fitting}
\label{sec_fit_2D}

\begin{figure}[t!]
\centerline
{
  \includegraphics[width=0.4\textwidth]{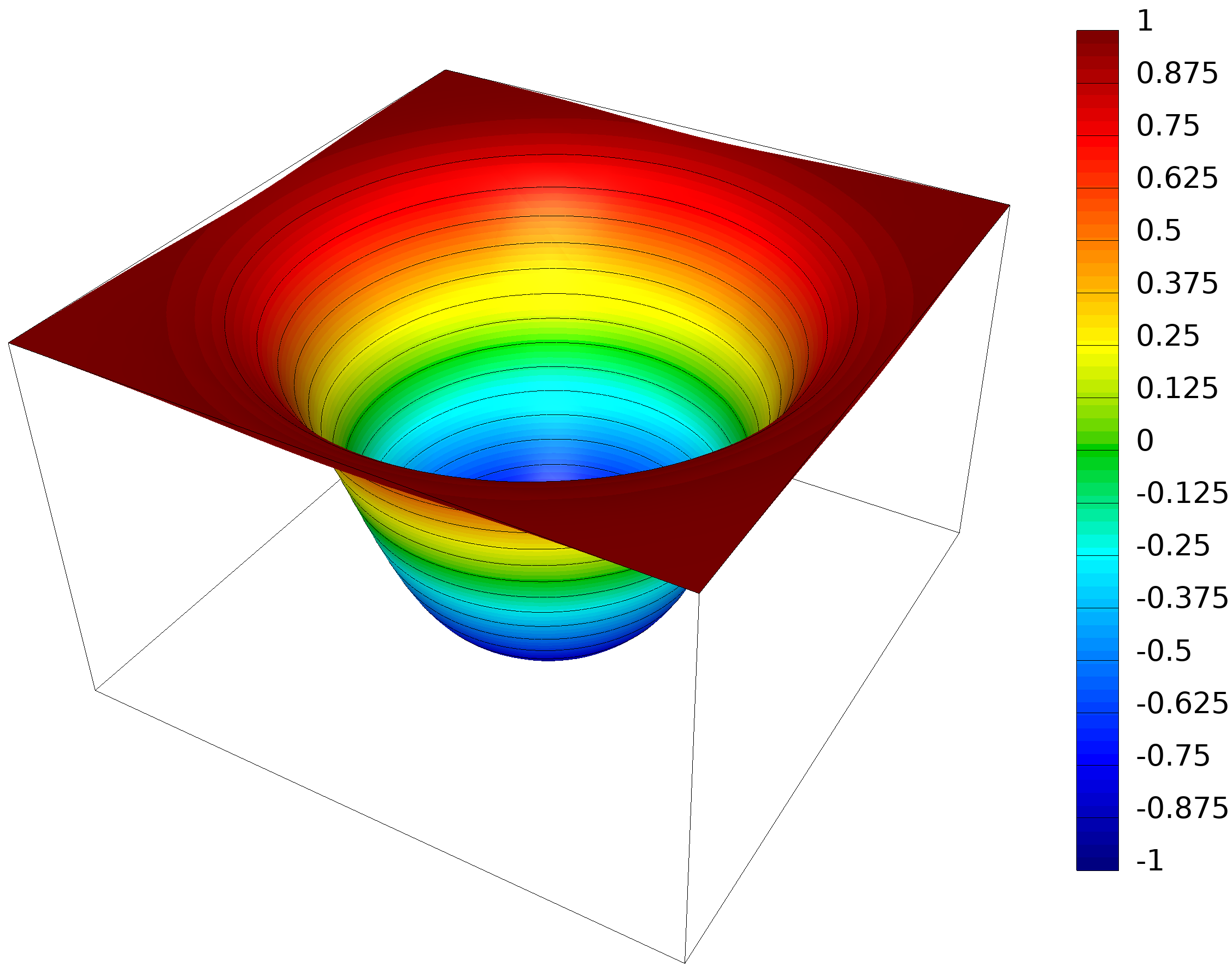}
}
\caption{Level set function $\sigma$ that is used for the
         2D surface fitting tests.}
\label{fig_2D_LS}
\end{figure}

\begin{figure}[t!]
\begin{center}
$\begin{array}{ccc}
\includegraphics[height=0.2\textwidth]{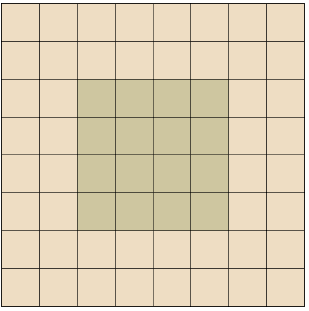} &
\includegraphics[height=0.2\textwidth]{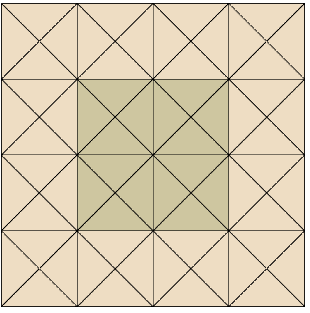}
\end{array}$
\end{center}
\vspace{-7mm}
\caption{Initial quadrilateral / triangle meshes and internal interfaces
         for the 2D surface fitting tests.}
\label{fig_2D_mesh}
\end{figure}

As a first example, we perform surface fitting to 3rd order 2D meshes consisting
of quads and triangles, see Figure \ref{fig_2D_mesh}.
For both meshes we use TMOP with tangential
relaxation and fitting capability \eqref{eq_F}-\eqref{eq_F_sigma} to adapt the
internal interface to the zero level-set of a function.
The zero level-set of the function $\sigma$ is defined such that it is located
at a distance of 0.3 from the center of the domain, $\bx_c = (0.5, 0.5)$,
see Figure \ref{fig_2D_LS}.
Although this level set is known analytically, the presented computations
represent and use $\sigma$ as a discrete finite element function.

\begin{figure}[b!]
\begin{center}
$\begin{array}{cc}
\includegraphics[height=0.2\textwidth]{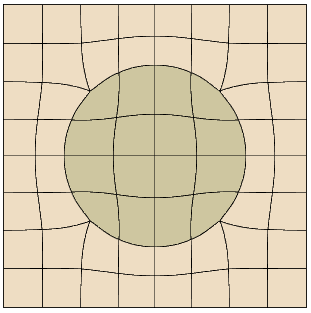} &
\includegraphics[height=0.2\textwidth]{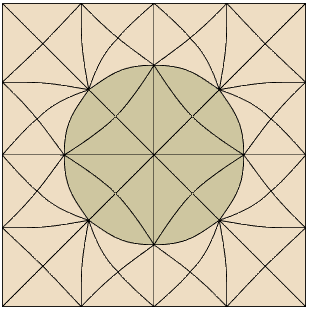} \\
\end{array}$
\end{center}
\vspace{-7mm}
\caption{Optimized quadrilateral and triangle meshes for the
         2D surface fitting tests.}
\label{fig_2D_mesh_opt}
\end{figure}

For both meshes we set $w_{\sigma} = 1000$.
The target matrix $T$ is chosen to represent an ideal element
(square or equilateral triangle),
and we use the shape-only metric
\[
\mu_{58}(T) = \frac{|T^tT|^2}{\tau^2}-{2|T|^2}{\tau}+2.
\]
The optimized meshes are shown in Figure \ref{fig_2D_mesh_opt}.
For the quadrilateral case, the TMOP objective function $F(\bx)$ decreases by $61.8\%$;
for the triangular case, $F(\bx)$ decreases by $71.9\%$.
For each of the optimized meshes, we also measure the average error for the
degree-of-freedoms associated with the material interface as
\begin{eqnarray}
\label{eq_zero_ls}
  e_{\mathcal{S}} = \frac{\sum_{s \in \mathcal{S}}
                    \bigg(||(\bxs-\bx_c)||_{2}-0.3\bigg)^2}{\sum_{s \in \mathcal{S}} 1},
\end{eqnarray}
where $\mathcal{S}$ is the set of nodes chosen to align to the interface,
$\bxs$ denotes the physical-location of the node $s$,
and $||\cdot||_{L_2}$ denotes the $L_2$-norm of a vector.
We observe that the average error at the zero-level set is
$O(10^{-6})$ for each of the meshes.

%-------

\subsection{3D Surface Fitting}
\label{sec_fit_3D}

\begin{figure}[b!]
\begin{center}
$\begin{array}{cc}
\includegraphics[height=0.23\textwidth]{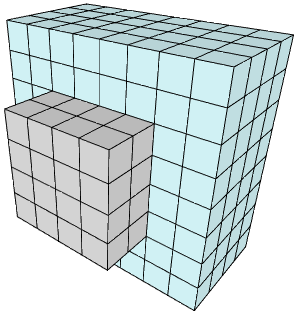} &
\includegraphics[height=0.23\textwidth]{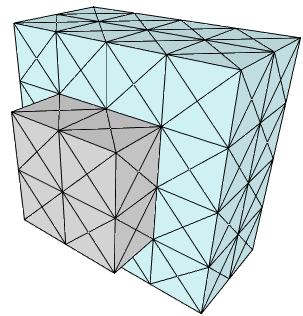} \\
\end{array}$
\end{center}
\vspace{-7mm}
\caption{Initial hexahedral / tetrahedral meshes and internal interfaces
         for the 3D surface fitting tests.}
\label{fig_3D_mesh}
\end{figure}

Next we consider the analogous 3D case and perform surface fitting to 2nd order
3D meshes consisting of hexes and tets.
Figure \ref{fig_3D_mesh} shows a sliced-view of the meshes, where the element
colors denote the two sides of the internal interface.
For both meshes the zero level set of the function $\sigma$ is
at a distance of 0.3 from the center of the domain $\bx_c = (0.5, 0.5, 0.5)$.
We set $w_{\sigma}=1000$.
The target matrix $T$ is chosen to represent an ideal element
(cube or equilateral tetrahedron) with size target
based on the size of the elements in the original mesh.
The mesh quality is controlled by the polyconvex shape+size metric $\mu_{333}(T)$.

\begin{figure}[b!]
\begin{center}
$\begin{array}{cc}
\includegraphics[height=0.22\textwidth]{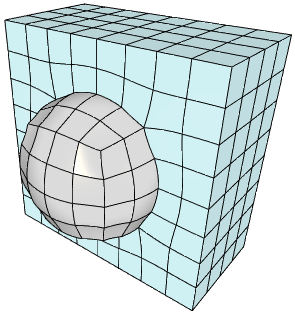} &
\includegraphics[height=0.22\textwidth]{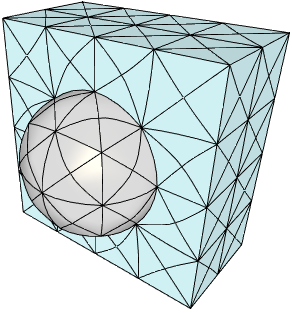} \\
\end{array}$
\end{center}
\vspace{-7mm}
\caption{Optimized hexahedral and tetrahedral meshes for the
         3D surface fitting tests.}
\label{fig_3D_mesh_opt}
\end{figure}

The optimized meshes are shown in Fig. \ref{fig_3D_mesh_opt}.
For the mesh with hexahedral elements, we observe that $F(\bx)$ reduces by $91.4\%$,
and the error at the zero level-set \eqref{eq_zero_ls}
is $O(10^{-5})$. For the mesh with tetrahedrons, $F(\bx)$ reduces by
$78.7\%$, with $e_{\mathcal{S}}=O(10^{-6})$.
The error at the zero level-set is relatively higher for the mesh with hexahedral elements
due to the mesh topology. The hex mesh has elements that have more than one
face at the material interface, which would have to tangentially flatten
along the zero level-set to exactly satisfy the surface fitting term \eqref{eq_F_sigma}.
Since this would result in an inverted mesh, the mesh optimizer fits the
mesh the best it can while ensuring mesh quality based on the
metric $\mu_{333}(T)$ used in the first term of \eqref{eq_F}.

%-------

\subsection{Taylor-Green Interface}
\label{sec_tg}

\begin{figure}[t!]
\centerline
{
  \includegraphics[width=0.3\textwidth]{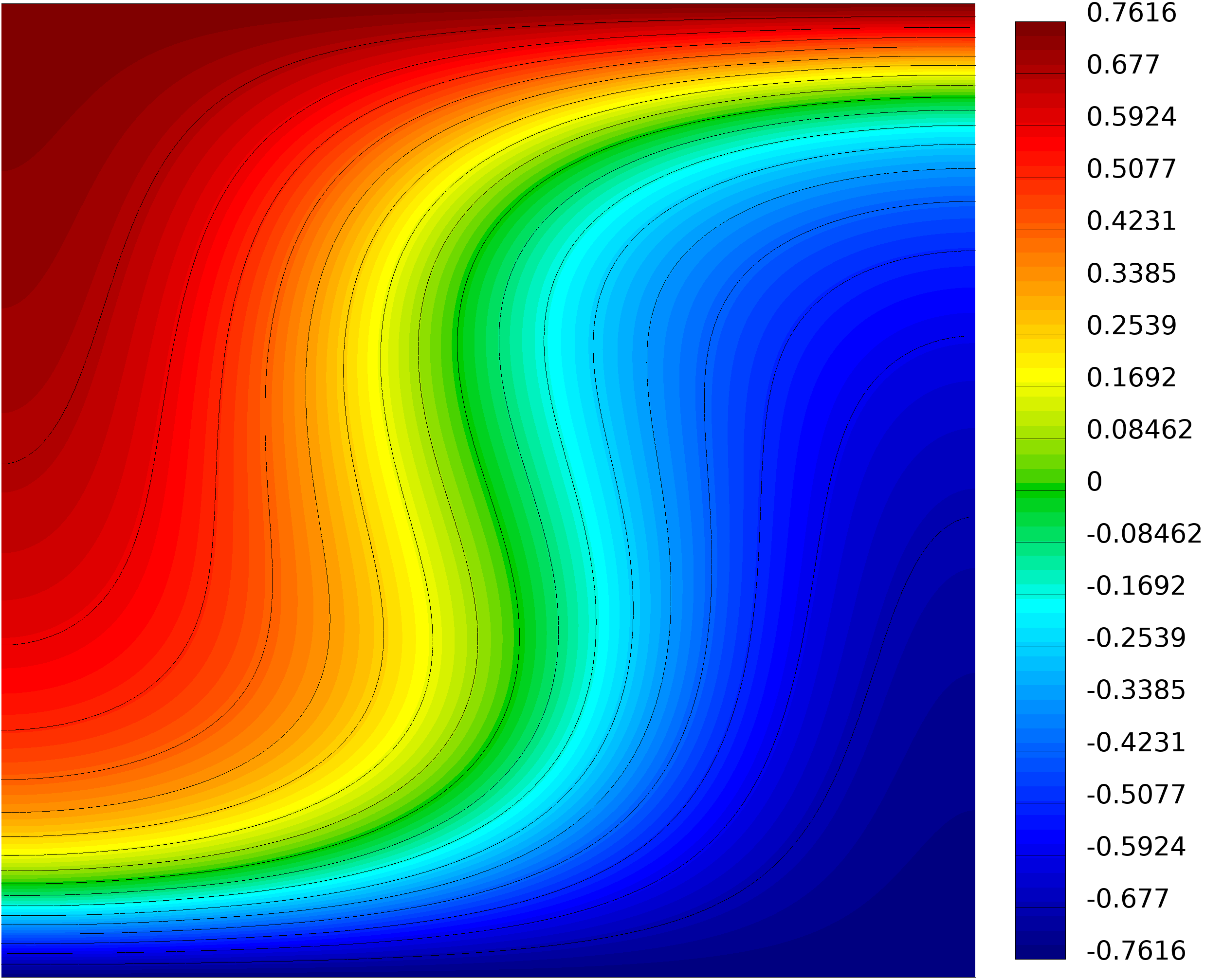}
}
\caption{Level sets of the finite element function $\sigma(\bx)$ for the
         material interface arising in the Taylor-Green test.}
\label{fig_tg_2d_ls}
\end{figure}

\begin{figure}[t!]
\centerline
{
  \includegraphics[width=0.22\textwidth]{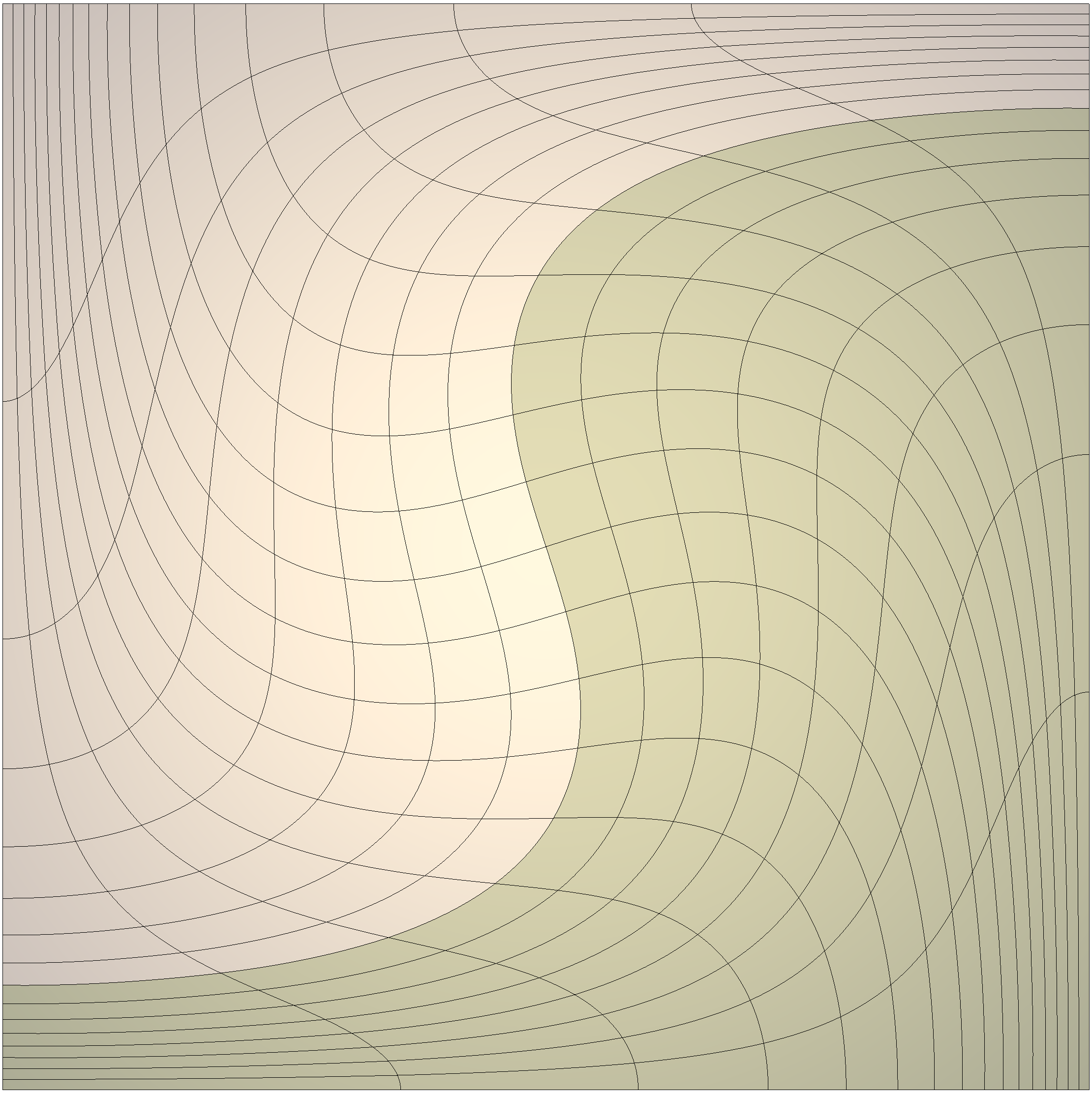}\hfil
  \includegraphics[width=0.22\textwidth]{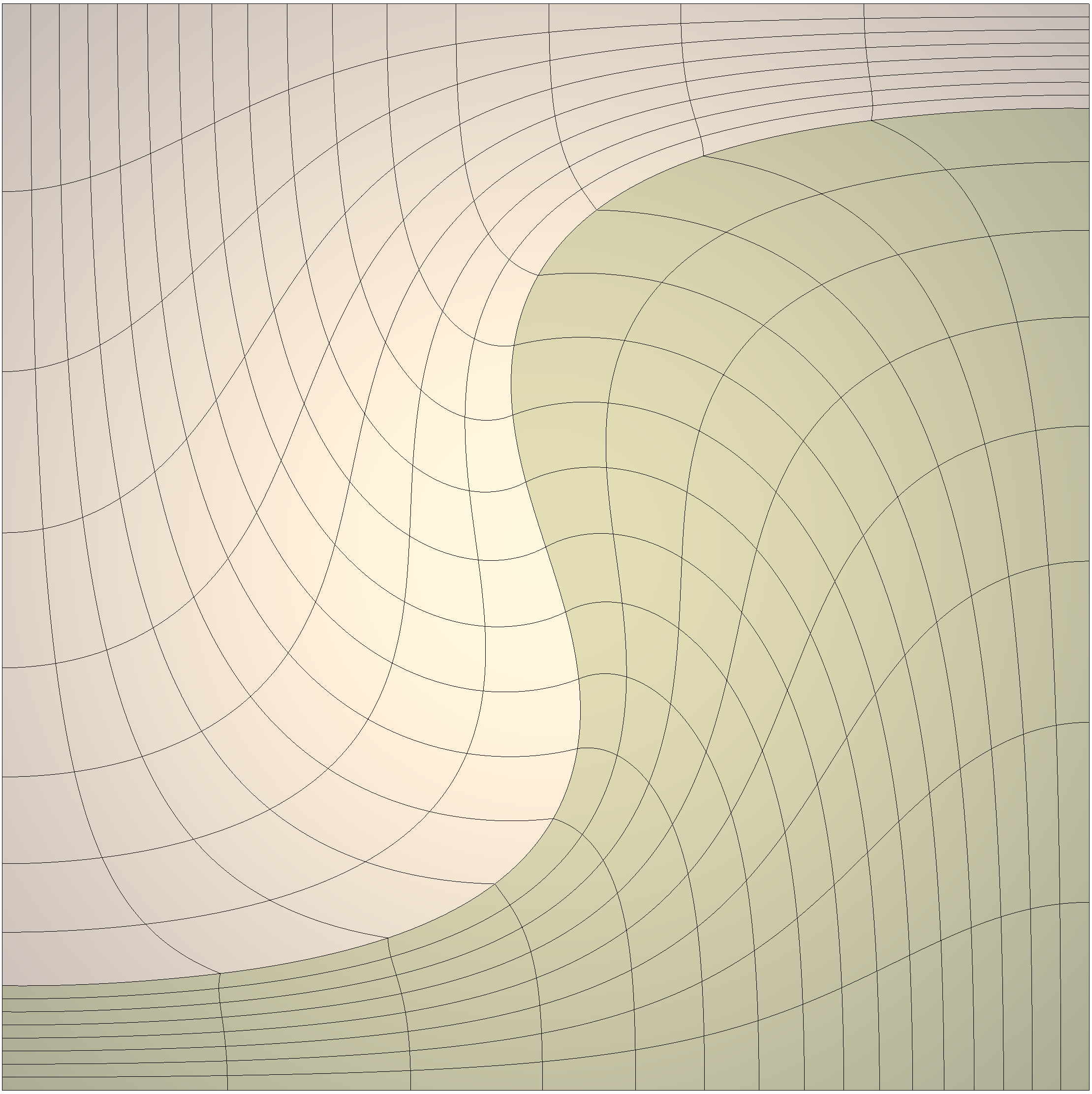}
}
\centerline
{
  \includegraphics[width=0.22\textwidth]{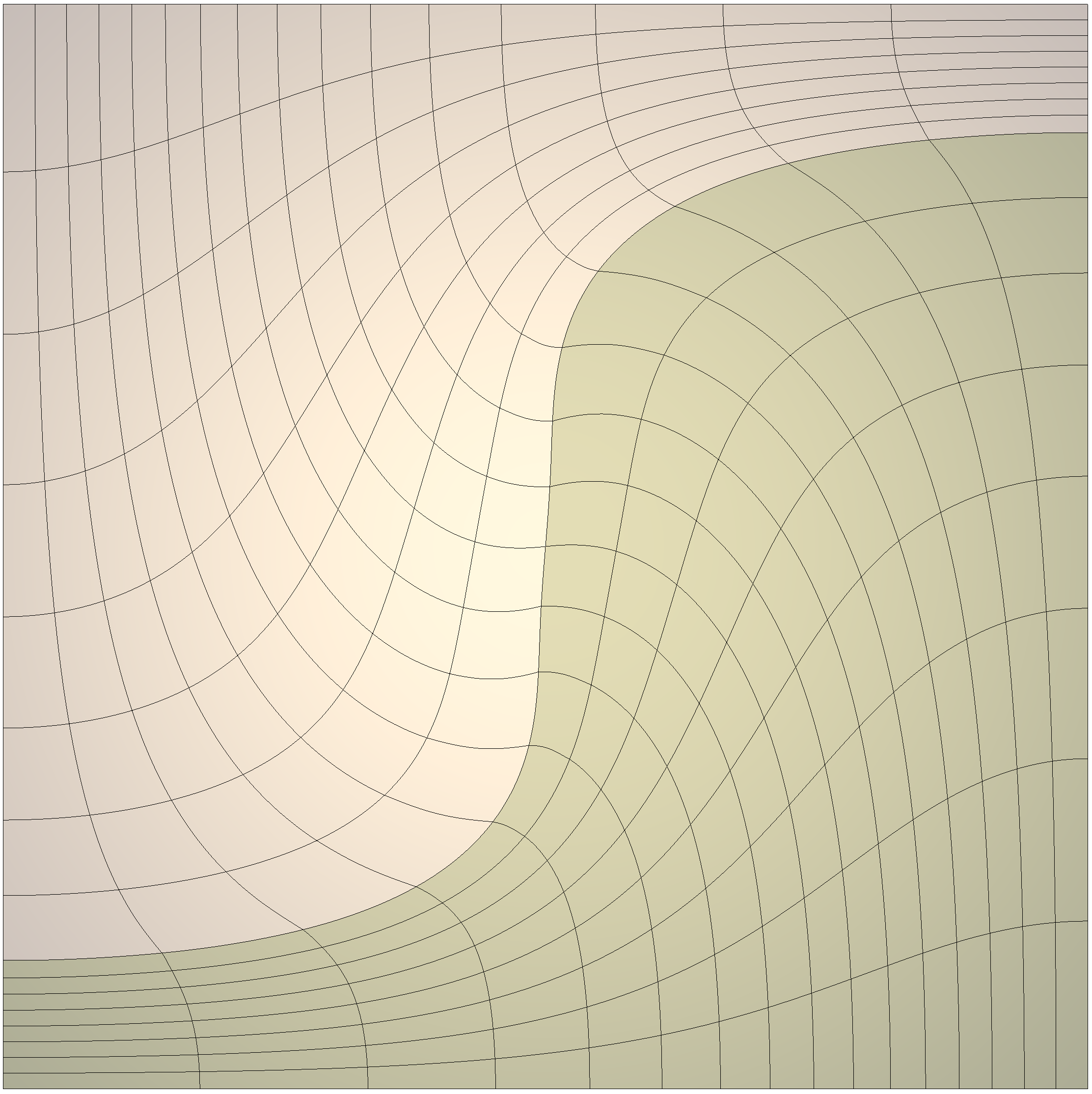}\hfil
  \includegraphics[width=0.22\textwidth]{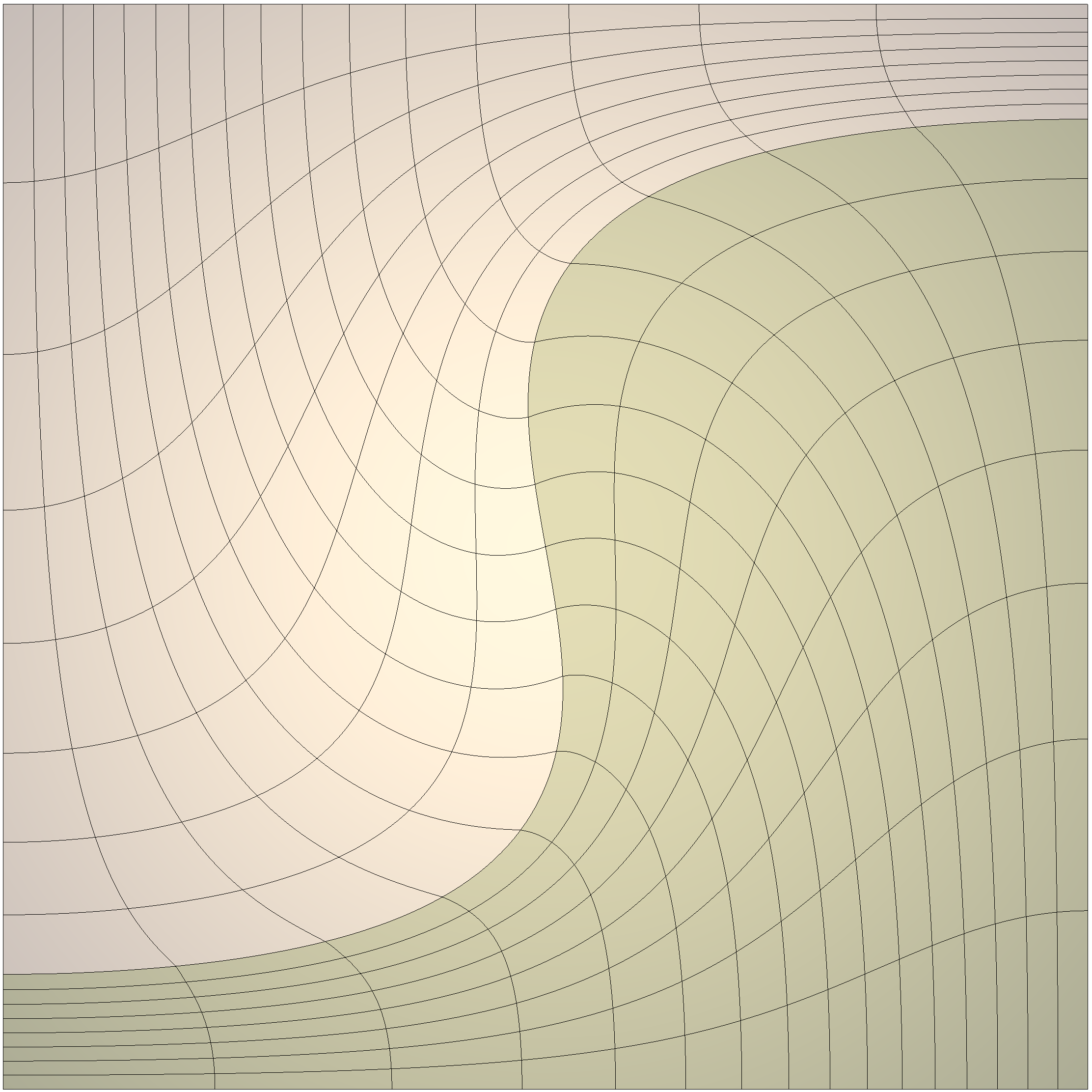}
}
\caption{Top row: initial mesh (left) and optimized mesh
         with fully constrained interface nodes (right).
         Bottom row: optimized mesh with $w_{\sigma} = 250$ (left) and
         optimized mesh with $w_{\sigma} = 1000$ (right).}
\label{fig_tg_2d}
\end{figure}

The primary target application of our work are
multimaterial moving mesh simulations.
In this context, a common goal of the mesh optimization procedure is to
improve the mesh while preserving an initially aligned material interface.
In this example we start with a 3rd order mesh and a material interface obtained
from a Taylor-Green moving mesh simulation \cite{BLAST2018}.
The level set function $\sigma(\bx)$ is shown in Figure \ref{fig_tg_2d_ls};
the initial mesh is shown in the top left panel of Figure \ref{fig_tg_2d}.

We compare results from three different setups, namely
(i) fixing the interface nodes during the optimization,
(ii) tangential relaxation with $w_{\sigma}=250$, and
(iii) tangential relaxation with $w_{\sigma}=1000$.
The decrease in the objective function and the interface position errors are
listed in Table \ref{tab_tg_2d}, while the optimized meshes for the three cases
are shown in Figure \ref{fig_tg_2d}.
In all cases the mesh quality is controlled by the quality metric $\mu_{80}$,
and the target Jacobians represent ideally shaped elements that maintain their
initial local size.
Since $\sigma$ is not known analytically, the interface position errors are
defined with respect to the violation of the zero level set:
\[
  \mathcal{E}_{avg} = \frac{1}{|\mathcal{S}|}
                      \sum_{s \in \mathcal{S}} \bar{\sigma}(\bxs),
  \quad
  \mathcal{E}_{max} = \max_{s \in \mathcal{S}} \bar{\sigma}(\bxs).
\]

\begin{table}[h!]
\begin{center}
  \begin{tabular}{c | c c c}
  \hline
      Approach & $F$ decrease & $\mathcal{E}_{avg}$ & $\mathcal{E}_{max}$ \\
  \hline
      Fixed interface   & 34.4\% & 0     & 0        \\
      $w_{\sigma}=250$  & 51.4\% & 8.2e-2 & 1.3e-1  \\
      $w_{\sigma}=1000$ & 42.6\% & 3.6e-2 & 6.4e-2  \\
  \hline
  \end{tabular}
\end{center}
\caption{Comparison of optimization strategies for the 2D Taylor-Green test.}
\label{tab_tg_2d}
\end{table}

We observe that the parameter $w_{\sigma}$ has the expected influence on the
fitting accuracy, namely, increasing its value leads to better agreement
with the zero level set of $\sigma$.
On the other hand, smaller values of $w_{\sigma}$ give more freedom to the
interface nodes, leading to better mesh quality at the price of increased error
of the interface positions.
Finally, completely fixing the interface nodes gives the lowest quality mesh,
as the optimizer has no control over the interfacial nodes.

\begin{figure}[t!]
\centerline
{
  \includegraphics[width=0.3\textwidth]{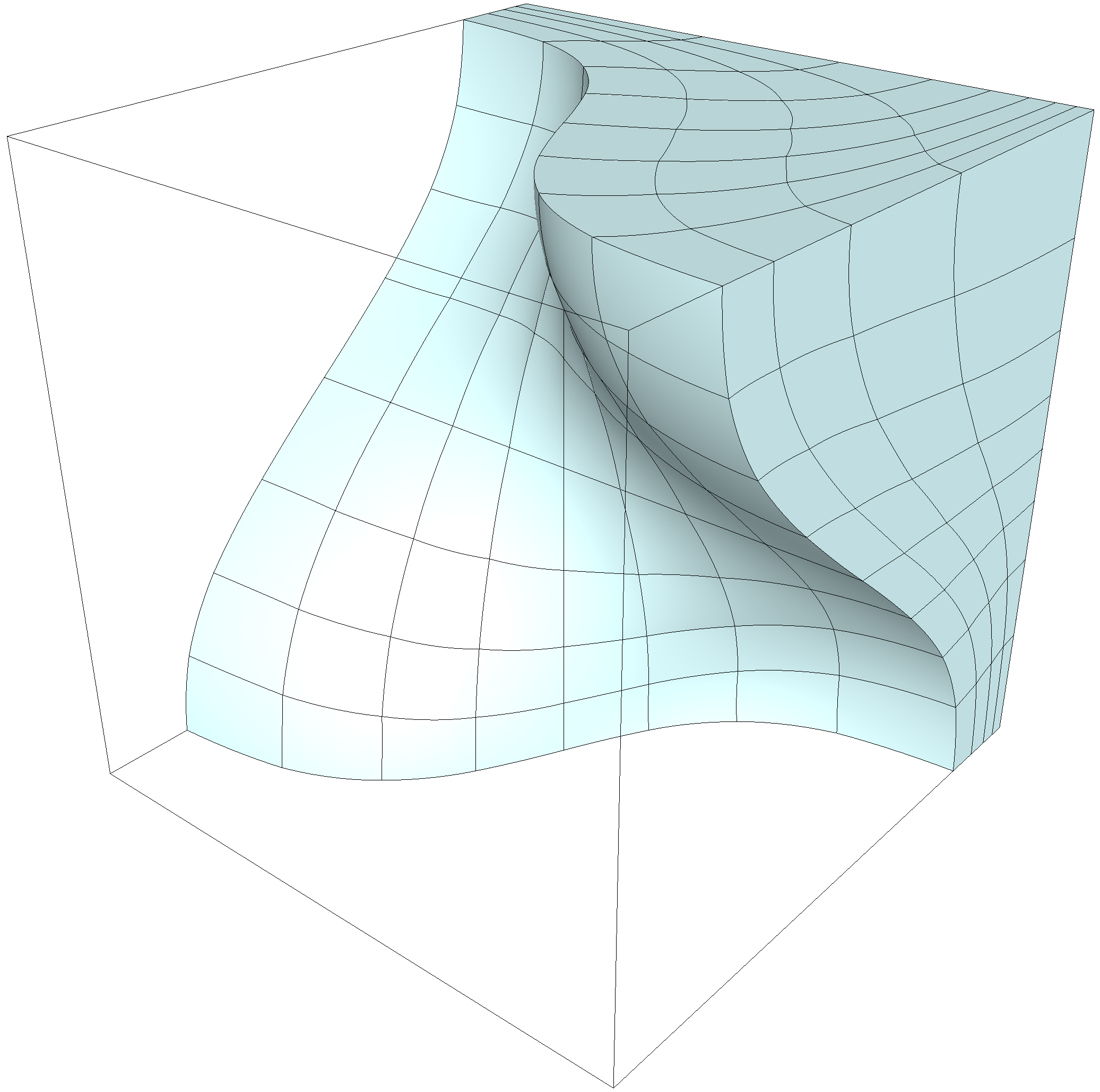}
}
\centerline
{
  \includegraphics[width=0.3\textwidth]{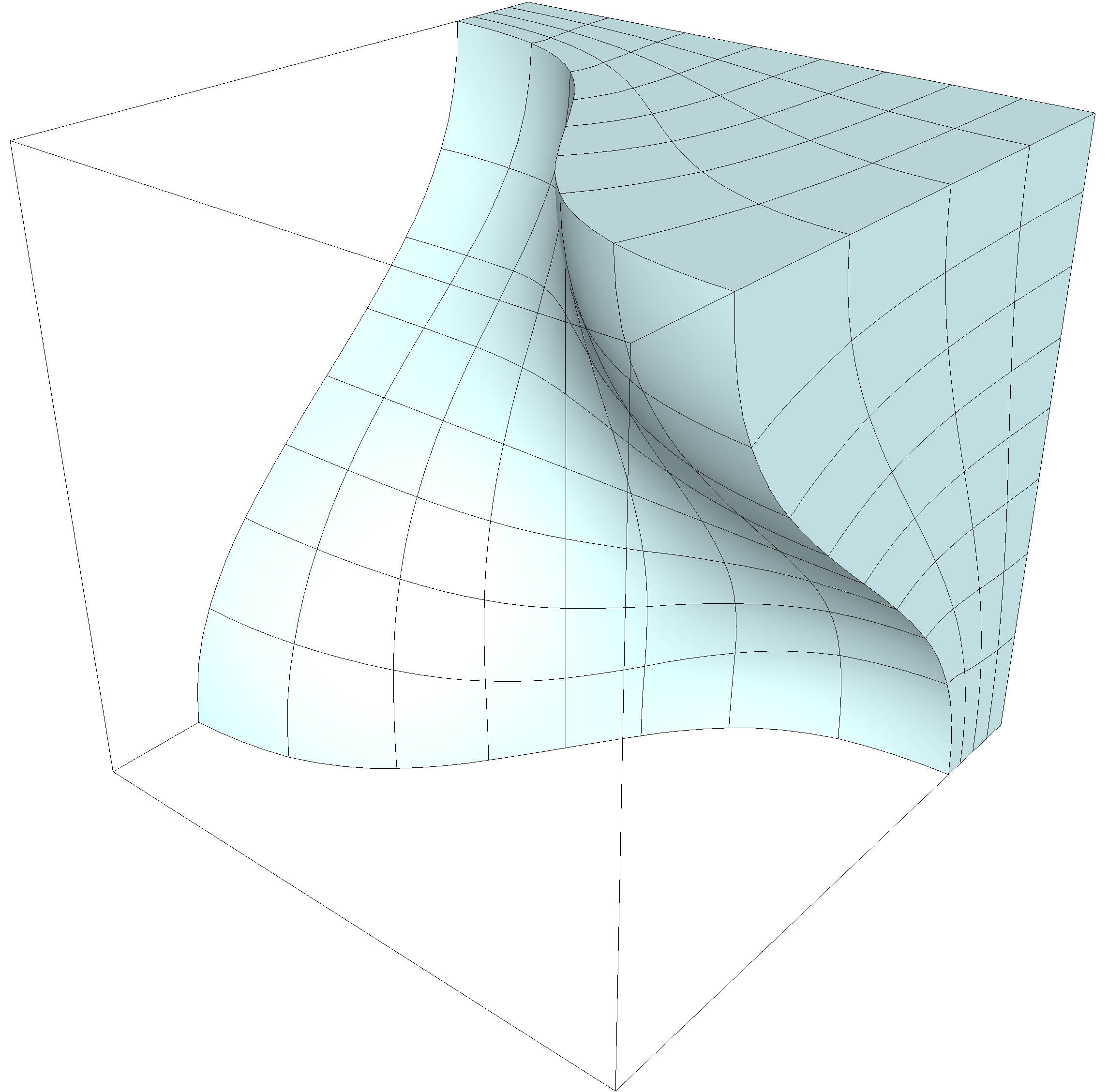}
}
\caption{Initial mesh (top) and optimized mesh with
         $w_{\sigma} = 1000$ (bottom) for the
         material interface arising in the 3D Taylor-Green test.}
\label{fig_tg_3d}
\end{figure}

The 3D version of the same problem is presented in Figure \ref{fig_tg_3d},
where we show cuts (along the material interface) of the initial mesh and the
optimized mesh with $w_{\sigma}=1000$.
The 3D comparison between the above three optimization strategies is listed in
Table \ref{tab_tg_3d}.
The 3D results confirm the observations made in the 2D tests.

\begin{table}[h!]
\begin{center}
  \begin{tabular}{c | c c c}
  \hline
      Approach & $F$ decrease & $\mathcal{E}_{avg}$ & $\mathcal{E}_{max}$ \\
  \hline
      Fixed interface   & 54.5\% & 0      & 0       \\
      $w_{\sigma}=250$  & 72.7\% & 4.9e-2 & 1.5e-1  \\
      $w_{\sigma}=1000$ & 63.6\% & 1.7e-2 & 6.2e-2  \\
  \hline
  \end{tabular}
\end{center}
\caption{Comparison of optimization strategies for the 3D Taylor-Green test.}
\label{tab_tg_3d}
\end{table}

%-------

\subsection{Rayleigh-Taylor Interface}
\label{sec_rt}

\begin{figure}[t!]
\centerline
{
  \includegraphics[width=0.092\textwidth]{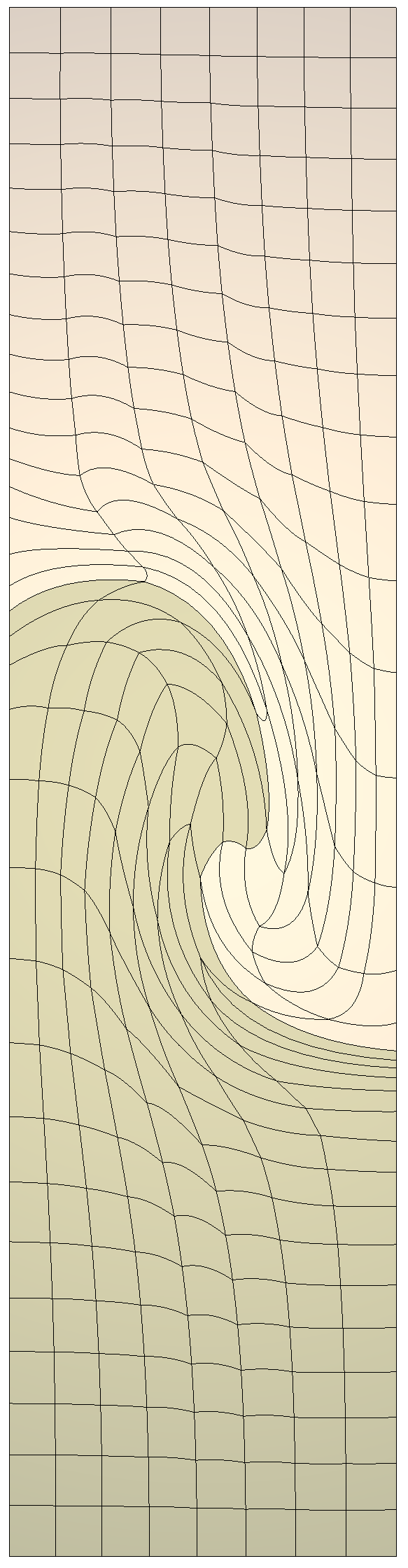}\hfil
  \includegraphics[width=0.09\textwidth]{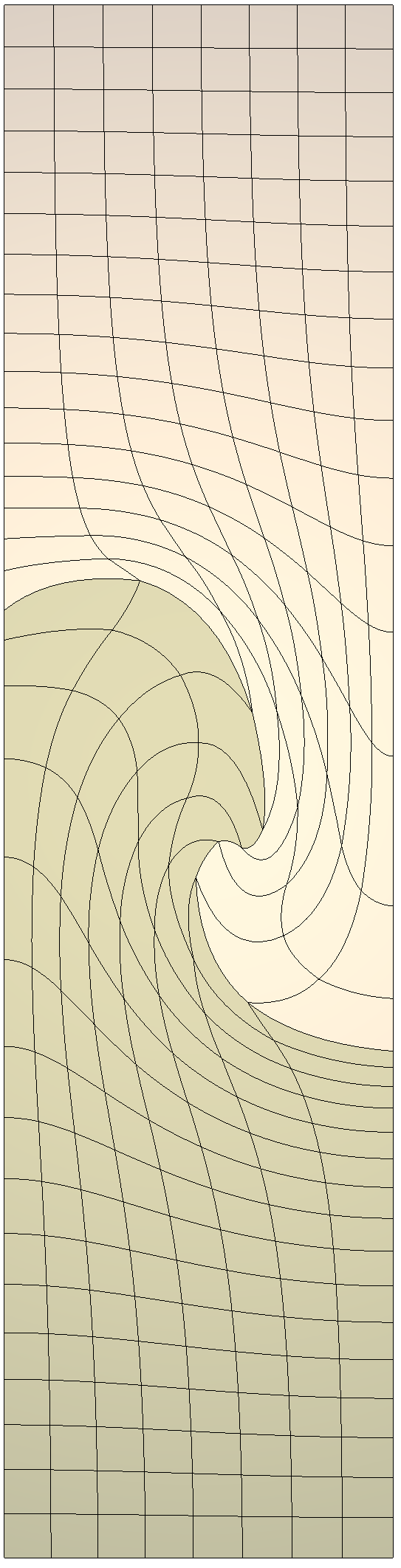}\hfil
  \includegraphics[width=0.09\textwidth]{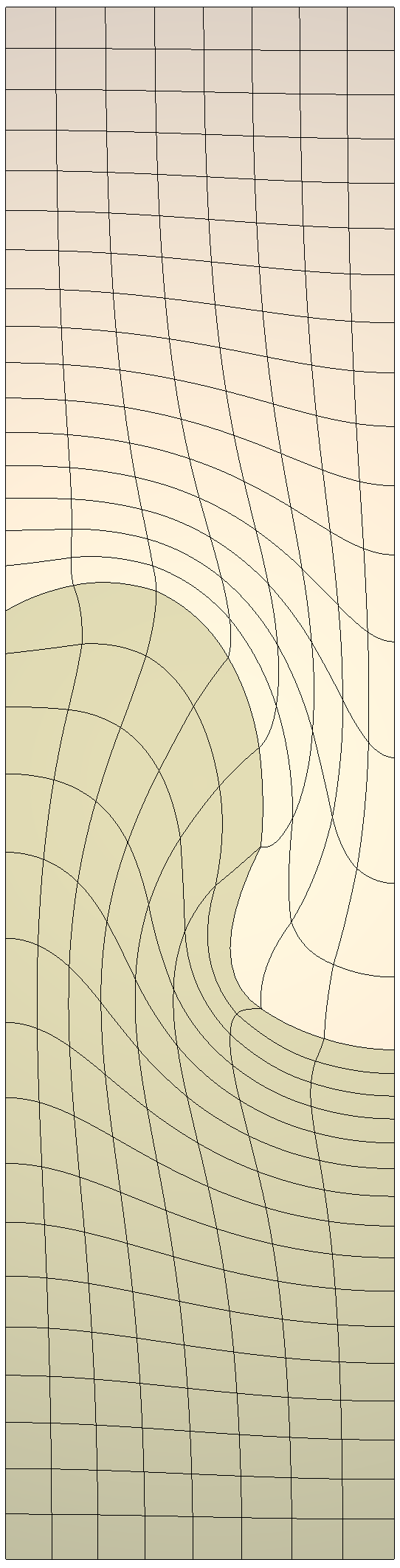}\hfil
  \includegraphics[width=0.153\textwidth]{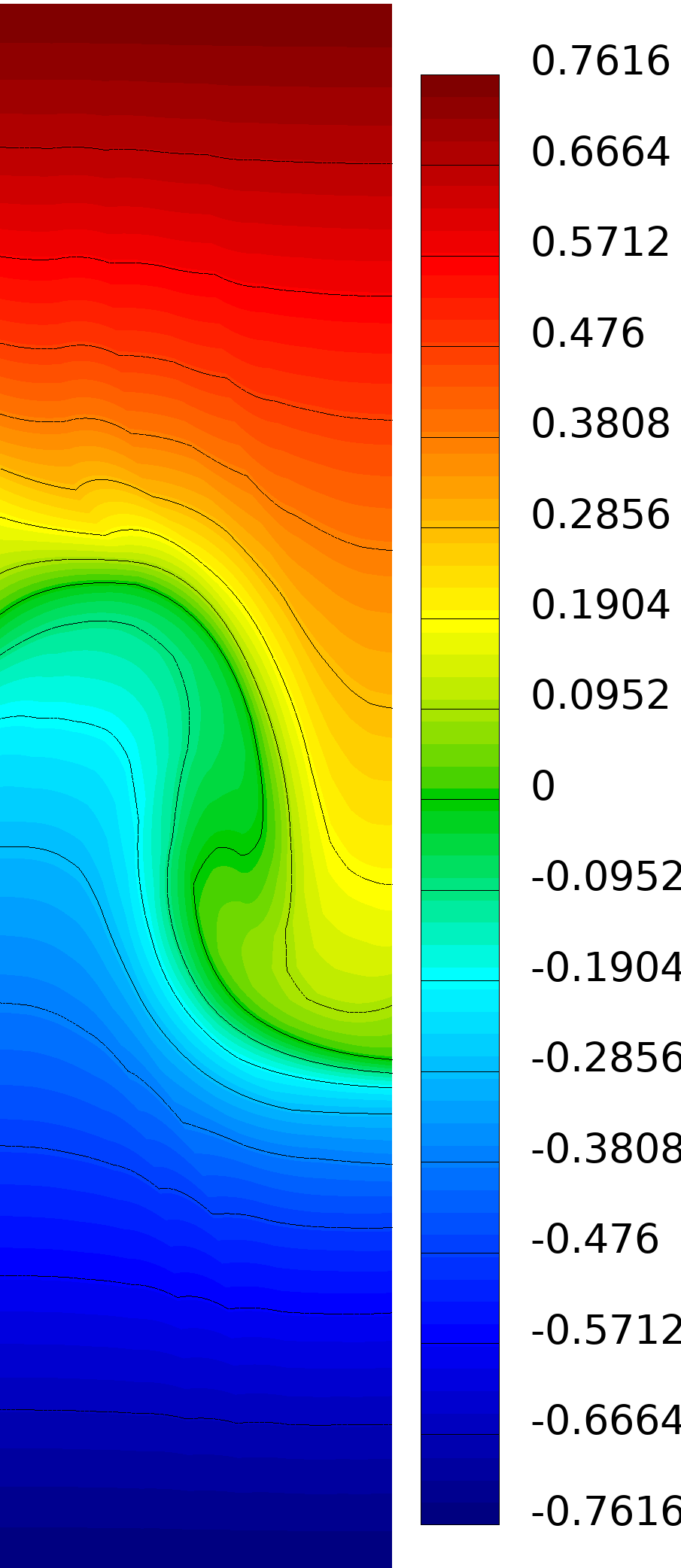}
}
\caption{Left to right: initial mesh,
         optimized mesh with fully constrained interface nodes,
         optimized mesh with $w_{\sigma}$ = 1.0e4, and
         level sets of the finite element function $\sigma(\bx)$
         for the 2D Rayleigh-Taylor test.}
\label{fig_rt_2d}
\end{figure}

Next we consider the Rayleigh-Taylor two-material problem \cite{Dobrev2012},
where the dynamics of the system leads to an interface that is not smooth.
The goal of this test is to demonstrate the behavior of the method for
non-smooth surfaces that contain fine local features.
The initial 2nd order mesh and the level set function
$\sigma(\bx)$ are shown in Figure \ref{fig_rt_2d}.

Again we utilize the $\mu_{80}$ quality metric and the target Jacobians
represent ideally shaped elements that maintain their initial local size.
Optimizing the mesh by exactly preserving the positions of the interfacial nodes
leads to a 38.6\% decrease in the objective function $F(\bx)$, while optimizing
with $w_{\sigma}$ = 1.0e4 leads to a decrease of 52.5\%, with
$\mathcal{E}_{avg}$ = 9.9e-3 and $\mathcal{E}_{max}$ = 4.4e-2.
The final optimized meshes for the two cases are shown in Figure \ref{fig_rt_2d}.

\begin{figure}[b!]
\centerline
{
  \includegraphics[width=0.14\textwidth]{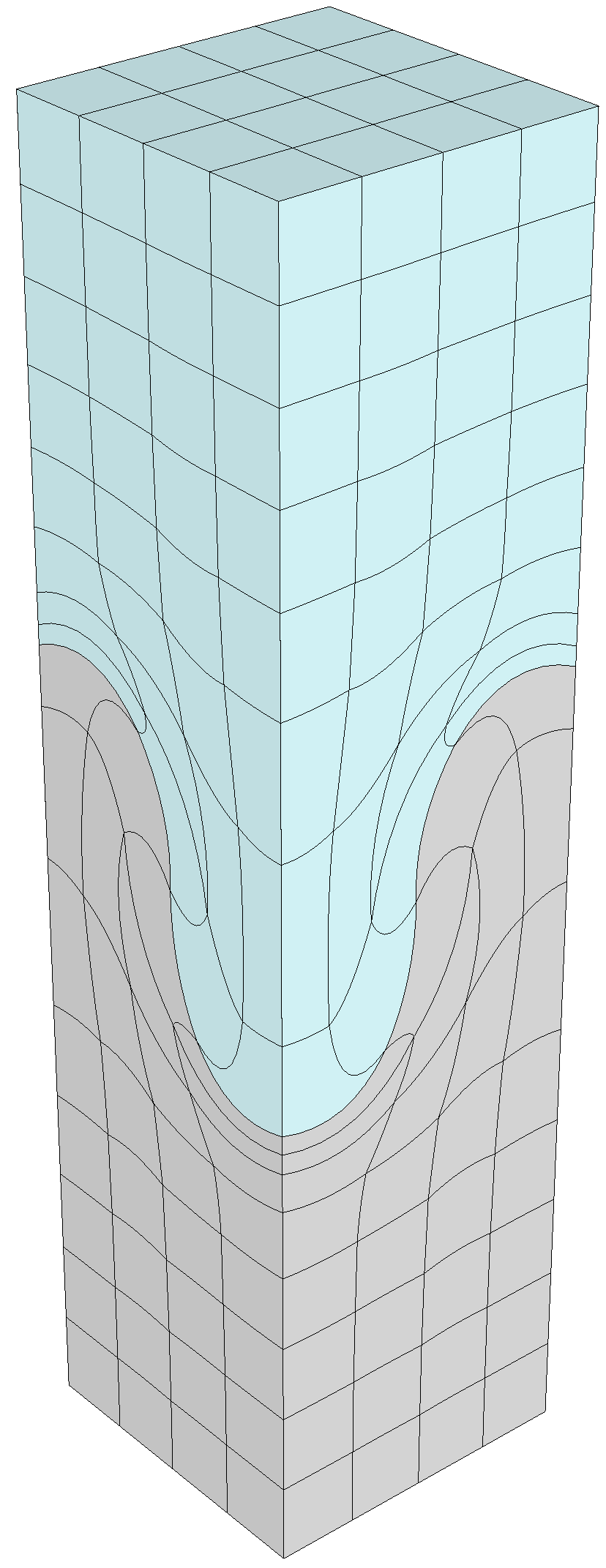}\hfil
  \includegraphics[width=0.24\textwidth]{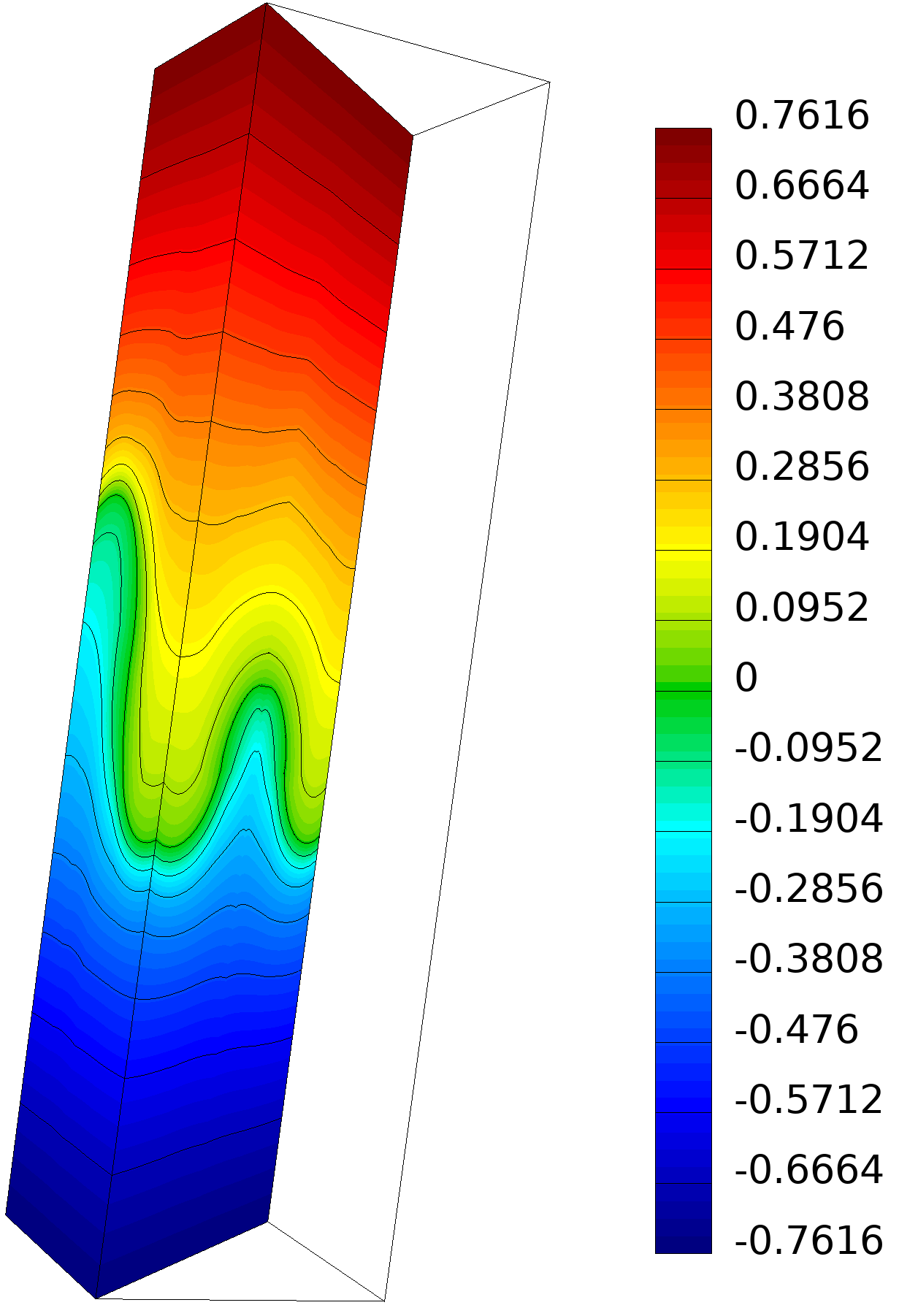}
  %% TK: the perspective in the right plot is confusing, can we use a straight version?
}
\caption{Initial mesh and cut of the level set function $\sigma(\bx)$
         for the 3D Rayleigh-Taylor test.}
\label{fig_rt_3d}
\end{figure}

The above results demonstrate that the proposed method can lead to
\textit{diffusion} of small surface features, even when $w_{\sigma}$ is large.
The reason is that the interface nodes $s \in \mathcal{S}$ are allowed to move
anywhere where the function $\sigma(\bxs)$ is small, as long as they improve
the mesh quality $\mu(\bx)$ in their local neighborhood.
Thus, small interface kinks as the ones visible in the left panel of
Figure \ref{fig_rt_2d} can be
eliminated without any significant resulting penalization in $F_{\sigma}$.
This behavior can be mitigated by detecting the local features of the surface,
e.g., by examining the gradients of $\sigma$, and performing localized treatments.
For example, one can fix only specific nodes or increase the resolution around
the feature by adapting the local size through $r$- or $h$-adaptivity \cite{dobrev2021hr}.
Another approach to address the non-smoothness of $\sigma$ is to replace the
function by $N_c$ smooth component functions $\sigma_1 \dots \sigma_{N_c}$,
and have a separate objective term for each smooth component in \eqref{eq_F}.
This strategy can be seen as the discrete analogy of the virtual geometry
surface decomposition models \cite{Roca2019}.
Exploring the above strategies will be the subject of future work.

\begin{figure}[t!]
\centerline
{
  \includegraphics[width=0.3\textwidth]{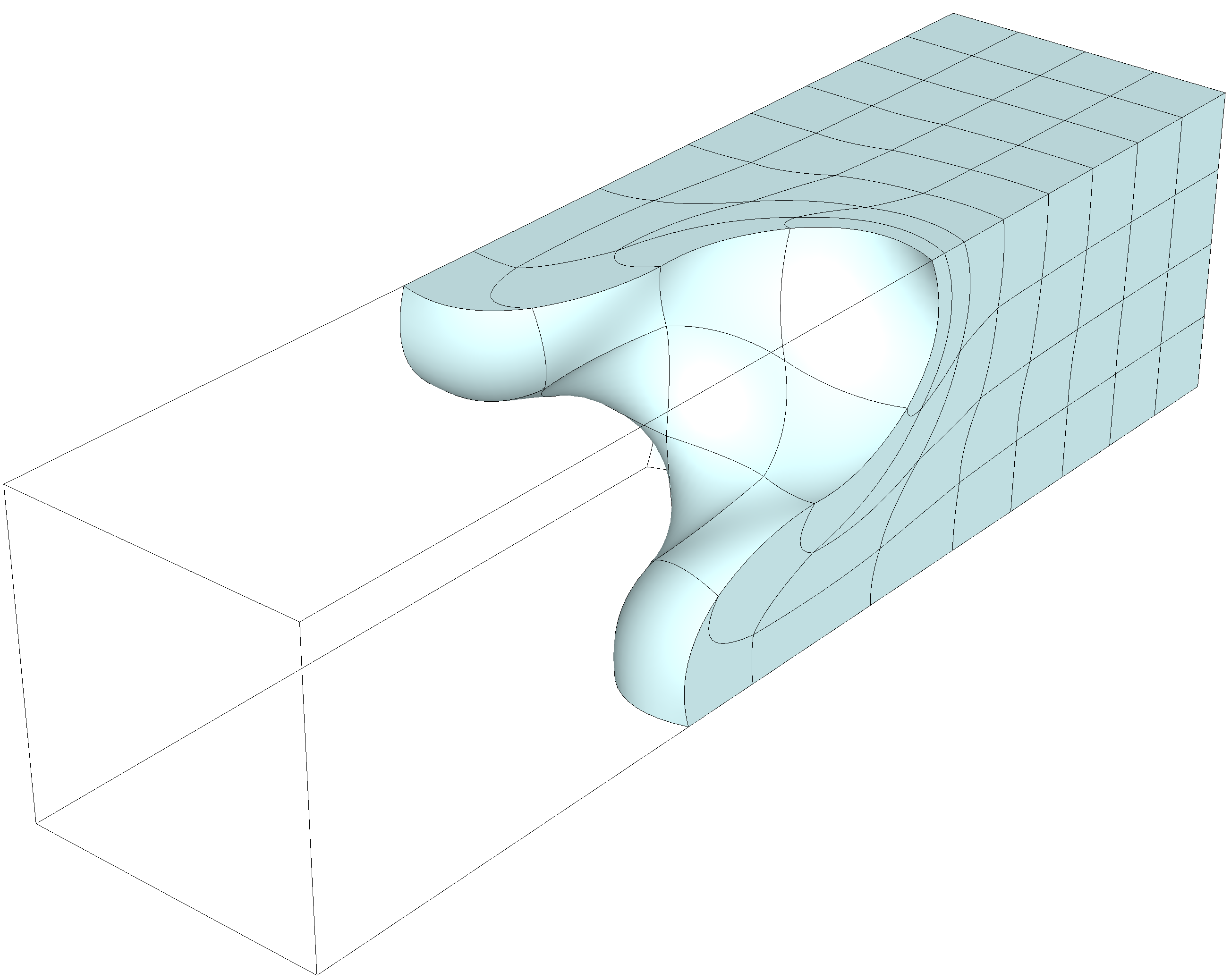}
}
\centerline
{
  \includegraphics[width=0.3\textwidth]{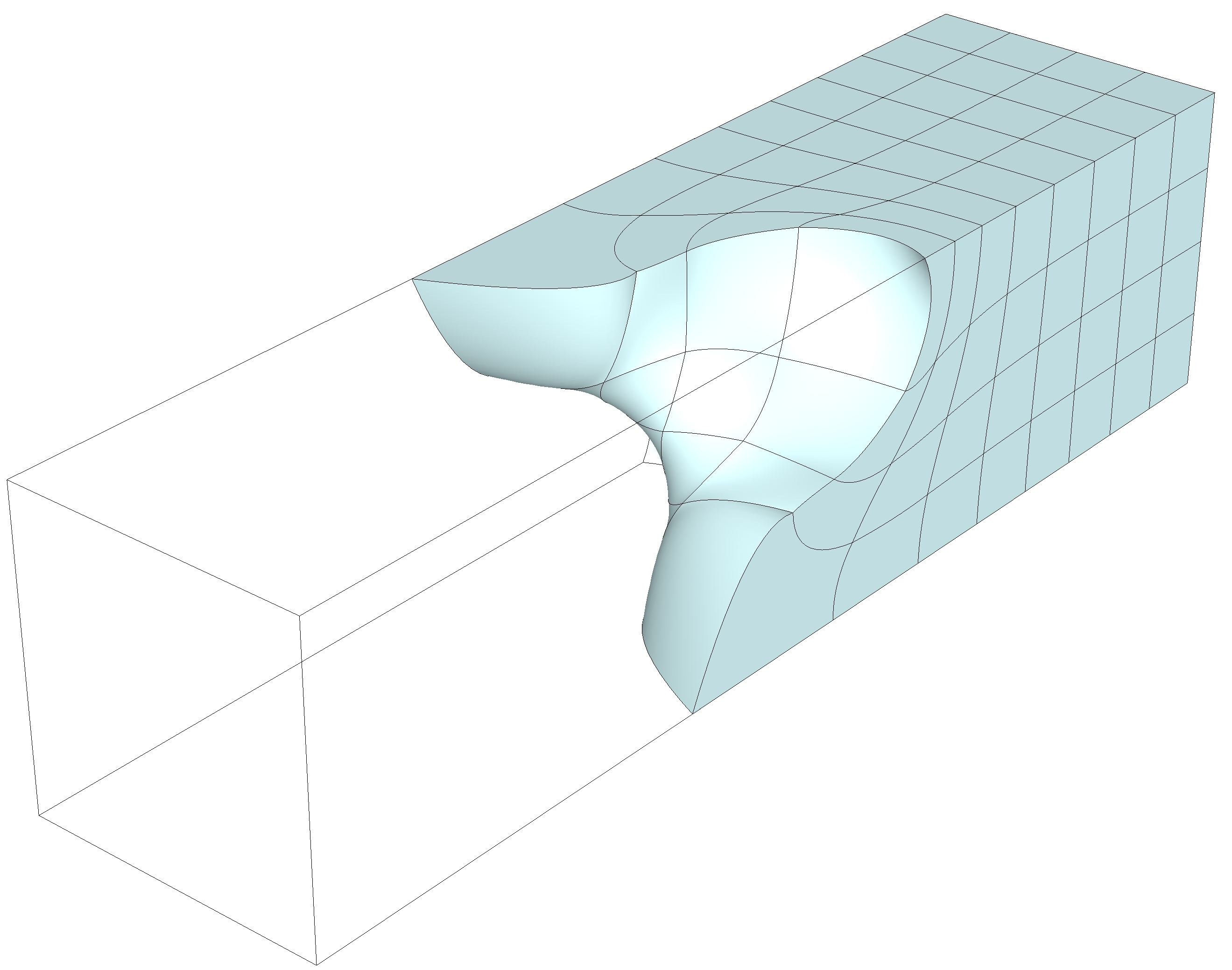}
}
\caption{Optimized mesh with fully constrained interface nodes (top) and
         optimized mesh with $w_{\sigma}$ = 1.0e4 (bottom) for the
         material interface arising in the 3D Rayleigh-Taylor test.}
\label{fig_tg_3d_cut}
\end{figure}

The 3D version of the same problem is presented in Figures \ref{fig_rt_3d} and
\ref{fig_tg_3d_cut}, where we observe behavior similar to the 2D case.
Optimizing the mesh by exactly preserving the positions of the interfacial nodes
leads to a 51.5\% decrease in the objective function $F(\bx)$, while optimizing
with $w_{\sigma}$ = 1.0e4 leads to a decrease of 62.9\%, with
$\mathcal{E}_{avg}$ = 8.8e-3 and $\mathcal{E}_{max}$ = 3.8e-2.

%-------------------------------------------------

\section{Conclusion}
\label{sec_concl}

We have presented a new approach to fit or align certain high-order mesh faces
to a surface given by a discrete level set function.
The alignment is imposed weakly by including a variational
penalty term in the objective functional.
The main advantage of the method is that its major steps can be implemented
strictly through finite element operations without performing any geometric
calculations, making the algorithm independent of dimension, mesh order,
and element type.
The main disadvantage of the algorithm is that local non-smooth surfaces
may be diffused in the optimization process.
Improving this aspect of the method will be future work, as discussed in
Section \ref{sec_rt}.

In the future the method will also be extended to handle the case of boundary
fitting, which is not addressed in this paper.
This will likely be approached by utilizing an auxiliary background mesh.
Extensive evaluation of the proposed method on larger problems with
complicated geometries will also be performed, which will require
robust strategies for combining several level set functions and more
involved node marking algorithms.

\paragraph{Disclaimer}
This document was prepared as an account of work sponsored by an agency of the
United States  government. Neither the United States government nor Lawrence
Livermore National Security,  LLC, nor any of their employees makes any
warranty, expressed or implied, or assumes any legal  liability or
responsibility for the accuracy, completeness, or usefulness of any information,
apparatus, product, or process disclosed, or represents that its use would not
infringe privately owned rights. Reference herein to any specific commercial
product, process, or service by trade name, trademark, manufacturer, or
otherwise does not necessarily constitute or imply its endorsement,
recommendation, or favoring by the United States government or Lawrence
Livermore National Security, LLC. The views and opinions of authors expressed
herein do not necessarily state or reflect those of the United States government
or Lawrence Livermore National Security, LLC, and shall not be used for
advertising or product endorsement purposes.

\bibliography{imr2021}

\end{document}